\documentclass[12pt]{article}
\usepackage{amsmath,amsxtra,amssymb,latexsym, amscd,amsthm}

\advance\voffset-1truecm\relax \advance\hoffset-1truecm\relax

scaled \magstep2  
scaled \magstep1

\def\N{{\Bbb N}}

\def\C{{\Bbb C}}

\def\Z{{\Bbb Z}}

\numberwithin{equation}{section}

\newcommand {\Cal}{\mathcal}

\newtheorem{lemma}{Lemma}[section]
\newtheorem{theorem}[lemma]{Theorem}
\newtheorem{corollary}[lemma]{Corollary}

\newtheorem{proposition}[lemma]{Proposition}

\def\leq{\leqslant}

\begin{document}
\title{A UNIQUENESS THEOREM\\ FOR MEROMORPHIC MAPS\\ WITH
MOVING HYPERSURFACES}

\author{$\quad$Gerd Dethloff and Tran Van Tan}
\date{$\quad$}
\maketitle
\vspace{-0.5cm}
\begin{abstract}
\noindent In this paper, we establish  a uniqueness theorem for
algebraically nondegenerate meromorphic maps of $\C^m$ into $\C P^n$
and slowly moving hypersurfaces  $Q_j \subset\C P^n,$ $j=1,\dots,q$  in (weakly) general position, where $q$ depends effectively on $n$ and on the degrees $d_j$ of the hypersurfaces
$Q_j$.
\end{abstract}

\section{Introduction}
One of the most striking consequences of Nevanlinna's theory
 was his ``five values'' theorem, which says that if $f$ and
$g$ are non-constant meromorphic functions on $\C$ such that
$f^{-1}(a_i)=g^{-1}(a_i)$ for five  distinct points $a_i$ in the
extended complex plane, then $f=g$. This theorem is an example of
what is now known as   ``uniqueness theorem". In 1975, Fujimoto
generalized this result of Nevanlinna to the case of
meromorphic maps of $\C^m$ into $\C P^n.$ In the last years, many uniqueness theorems for meromorphic maps with
hyperplanes (both for fixed and for moving ones) have
been established. 

\hbox to 5cm {\hrulefill }

 \noindent {\small Mathematics Subject
Classification 2000: Primary 32H30; Secondary 32H04, 32H25,
14J70.}

\noindent {\small Key words:  Nevanlinna theory,
Second Main Theorem, Uniqueness Theorem.}

 {\small The first named author was partially supported by the Fields Institute Toronto. The second named author was partially supported by the
post-doctoral research program of the Abdus Salam International Centre for Theoretical Physics.}

For the case of hypersurfaces, however, there are
so far only the uniqueness theorem of Thai and Tan \cite{ThT} for the case
of Fermat moving hypersurfaces  and the one of Dulock and Ru
\cite{DR2} for the case of (general) fixed hypersurfaces. More precisely, in \cite{DR2}, 
Dulock and Ru prove that one has a uniqueness theorem for algebraically non-degenate
holomorphic maps $f,g:\C \rightarrow \C P^n$ 
satisfying  $f=g$ on $\cup_{i=1}^q(f^{-1}(Q_i)\cup g^{-1}(Q_i)),$
with respect to 
$q > (n+1)+\frac{2Mn}{\tilde d} + \frac{1}{2}$ fixed hypersurfaces  $Q_i \subset \C P^n$ in general position, where $\tilde d$ is the minimum of the degrees of these hypersurfaces and $M$ is the truncation level
in the Second Main Theorem for fixed hypersurface targets obtained by An-Phuong \cite{AP}
with $\epsilon = \frac{1}{2}$. Their method of proof comes from their paper \cite{DR1},
where they prove a uniqueness theorem for holomorphic curves into abelian varieties.

In this paper,  by a method different to the one used by Dulock and Ru, we prove a uniqueness theorem for the
case of slowly moving hypersurfaces (Corollary \ref{3.2} below). 
More precisely, we prove that one has a uniqueness theorem for algebraically non-degenate
meromorphic maps $f,g:\C^m \rightarrow \C P^n$ 
satisfying  $f=g$ on $\cup_{i=1}^q(f^{-1}(Q_i)\cup g^{-1}(Q_i))$ 
with respect to 
$q > (n+1) + \frac{2nL}{\tilde d} + \frac{1}{2}$ moving hypersurfaces $Q_i \subset \C P^n$
in (weakly) general position, where $\tilde d$ is the minimum of the degrees of these hypersurfaces and $L$ is the truncation level
in the Second Main Theorem for moving hypersurface targets obtained by the authors in   \cite{DT1} with $\epsilon = \frac{1}{2}$.
Moreover, under the additional assumption that the $f^{-1}(Q_i)$, $i=1,...,q$ intersect properly, $q > (n+1) + \frac{2L}{\tilde d} + \frac{1}{2}$ moving hypersurfaces are sufficient.
We remark that in the special case of
fixed hypersurfaces, our result gives back the uniqueness theorem
of Dulock and Ru (remark that $L \leq M$ in this case).
Moreover, we give our uniqueness
theorem in a slightly more general form (Theorem \ref{3.1} below), requiring assumptions on the $(p-1)$ first derivatives of the maps, which gives in return a better bounds 
 on the number of moving hypersurfaces in $\C P^n$, namely
 $q > (n+1) + \frac{2nL}{p\tilde d} + \frac{1}{2}$ respectively $q > (n+1) + \frac{2L}{p\tilde d} + \frac{1}{2}$.
 
\section{Preliminaries}
 For $z = (z_1,\dots,z_m) \in \C^m$, we set
$\Vert z \Vert = \Big(\sum\limits_{j=1}^m |z_j|^2\Big)^{1/2}$ and
define
\begin{align*}
B(r) &= \{ z \in \C^m : \Vert z \Vert < r\},\quad
S(r) = \{ z \in \C^m : \Vert z \Vert = r\},\\
d^c &= \dfrac{\sqrt{-1}}{4\pi}(\overline \partial - \partial),\quad
\Cal V = \big(dd^c \Vert z\Vert^2\big)^{m-1},\; \sigma = d^c
\text{log}\Vert z\Vert^2 \land \big(dd^c\text{log}\Vert z
\Vert\big)^{m-1}.
\end{align*}
Let $L$ be a positive integer or $+\infty$ and $\nu$ be a divisor on
$\C^m.$ Set $ |\nu| = \overline {\{z : \nu(z) \neq 0\}}.$
 We define the counting function of $\nu$ by
\begin{align*}
N^{(L)}_\nu(r) := \int\limits_1^r \frac{n^{(L)}(t)}{t^{2m-1}}dt\quad
(1 < r < +\infty),
\end{align*}
where
\begin{align*}
n^{(L)}(t) &= \int\limits_{|\nu | \cap B(t)} \text{min}\{\nu
,L\}\cdot \Cal V\ \quad
\text{for}\quad m \geq 2 \ \text{and}\\
n^{(L)}(t) &= \sum_{|z| \leq t}\text{min}\{ \nu(z),L\} \qquad\quad
\text{for}\quad m = 1.
\end{align*}

Let $F$ be a nonzero holomorphic function on $\C^m$. For a set
$\alpha = (\alpha_1,\dots,\alpha_m)$ of nonnegative integers, we set
$|\alpha| := \alpha_1 + \dots + \alpha_m$ and $\mathcal D^\alpha F
:= \dfrac{\partial^{|\alpha|}} {\partial^{\alpha_1}z_1 \cdots
\partial^{\alpha_m}z_m}\,\cdotp$
 We define the zero divisor  $\nu_F$ of $F$ by
\begin{align*}
\nu_F(z) = \max \big\{ p : \mathcal D^\alpha F(z) = 0 \ \text{for
all $\alpha$ with}\ |\alpha| < p \big\}.
\end{align*}

Let $\varphi$ be a nonzero meromorphic function on $\C^m$. The zero
divisor $\nu_\varphi$ of $\varphi$ is defined as follows: For each
$a \in \C^m$, we choose nonzero holomorphic functions $F$ and $G$ on
a neighborhood $U$ of $a$ such that $\varphi = \dfrac{F}{G}$ on $U$
and $\text{dim}\big(F^{-1}(0) \cap G^{-1}(0)\big) \leq m-2$, then we
put $\nu_\varphi(a) := \nu_F(a)$.

Set $N_\varphi^{(L)}(r):=N_{\nu_\varphi}^{(L)}(r).$ For brevity we
will omit the character ${}^{(L)}$ in the counting function if
$L=+\infty.$

Let $f$ be a meromorphic map of $\C^m$ into $\C P^n$. For arbitrary
fixed homogeneous coordinates $(w_0: \cdots : w_n)$ of $\C P^n$, we
take a reduced representation $f = (f_0 : \cdots : f_n)$, which
means that each $f_i$ is a holomorphic function on $\C^m$ and $f(z)
= (f_0(z) : \cdots : f_n(z))$ outside the analytic set $\{ z :
f_0(z) = \cdots = f_n(z) = 0\}$ of codimension $\geq 2$. Set $\Vert
f  \Vert = \max \{ |f_0|, \dots  , |f_n| \}$.

The characteristic function of $f$ is defined by
\begin{align*}
T_f(r) := \int\limits_{S(r)}\text{log}\Vert f \Vert \sigma -
\int\limits_{S(1)} \text{log}\Vert f \Vert \sigma ,\quad 1 < r < +\infty.
\end{align*}

 For a meromorphic function $\varphi$ on $\C^m$, the characteristic function
$T_\varphi(r)$ of $\varphi$ is defined by considering  $\varphi$ as
a meromorphic map of $\C^m$ into $\C P^1$.

Let $f$ be a nonconstant meromorphic map of $\C^m$ into $\C P^n$.
We say that a meromorphic function $\varphi$ on $\C^m$ is ``small"
with respect to $f$ if $T_\varphi(r) = o(T_f(r))$ as $r \to \infty$
(outside a set of finite Lebesgue measure).

Denote by $\mathcal M$ the field of all meromorphic functions on
$\C^m$ and by $\Cal K_f$ the subfield of $\mathcal M$  which consists of 
all ``small" (with respect to $f$) meromorphic functions on $\C^m$.

For a homogeneous polynomial $Q \in \Cal M[x_0,\dots,x_n]$ of degree
$d \geq 1$ 
we write $Q = \sum\limits_{I \in \Cal T_{d}} a_{I}x^I,$ where $ \Cal
T_d := \big\{ (i_0,\dots,i_n) \in \N_0^{n+1} : i_0 + \dots + i_n = d
\big\}$ and $x^I = x_0^{i_0} \cdots x_n^{i_n}$ for $x =
(x_0,\dots,x_n)$ and $I = (i_0, \dots,i_n) \in \Cal
T_d.$
Denote by $Q(z)= Q(z)(x_0, \dots, x_n)=\sum\limits_{I \in \Cal T_{d}} a_{I}(z)x^I$ the homogeneous  polynomial over $\C$ obtained by
evaluating the coefficients of $Q$ at a specific point $z \in \C^m$
in which all coefficient functions of $Q$ are holomorphic.

Let $Q \in \Cal M[x_0,\dots,x_n]$ of degree
$d \geq 1$ 
with $Q(f):=Q(f_0, \dots , f_n) \not\equiv 0$. We define
$$
N^{(L)}_f(r,Q) := N^{(L)}_{Q(f)}(r)\;\;\text{and}\;\;
f^{-1}(Q):=\{z: \nu_{Q(f)}>0\}.$$

The \textit{First Main Theorem} of Nevanlinna theory gives, for
$Q = \sum\limits_{I \in \Cal T_{d}} a_{I}x^I$ with 
$Q(f):=Q(f_0, \dots , f_n) \not\equiv 0$ :
%$$ N(r, Q)\leq d\cdot T_f(r)+ O\big(\sum_{I_1, I_2\in \Cal
%T_d}T_{\frac{a_{I_1}}{a_{I_2}}}(r)\big).$$
$$ N(r, Q)\leq d\cdot T_f(r)+ O\big(\sum_{I \in \Cal
T_d}T_{{a_{I}}}(r)\big).$$

Let
\begin{align*}
Q_j = \sum\limits_{I \in \Cal T_{d_j}} a_{jI}x^I \quad (j = 1,\dots,q)
\end{align*}
 be homogeneous polynomials in $\Cal K_f[x_0,\dots,x_n]$ with
$\text{deg}\,Q_j = d_j \geq 1.$  Denote by $\Cal
K_{\{Q_j\}_{j=1}^q}$ the field over $\C$ of   all meromorphic
functions on $\C^m$ generated by all quotients
$\big\{\frac{a_{jI_{1}}}{a_{jI_{2}}} :a_{jI_{2}}\not\equiv 0,
I_{1}, I_{2}
 \in \Cal T_{d_j}; j \in \{1,\dots,q\} \big\}$.
We say that $f$ is algebraically nondegenerate over $\Cal
K_{\{Q_j\}_{j=1}^q}$ if there is no nonzero homogeneous polynomial
$Q \in \Cal K_{\{Q_j\}_{j=1}^q}[x_0,\dots,x_n]$  such that
$Q(f_0,\dots,$ $f_n) \equiv 0$.

We say that a set $\{Q_j\}_{j=1}^q$ $(q \geq n+1)$ of homogeneous
polynomials in $\Cal K_f [x_0,\dots,$ $x_n]$ is admissible (or in (weakly) general position)
if there exists $z \in \C^m$
in which all coefficient functions of all $Q_j$, $j=1,...,q$ are holomorphic and such that for any
$1 \leq j_0 < \dots < j_n \leq q$ the system of equations
\begin{align} \label{zz}
\left\{ \begin{matrix}
Q_{j_i}(z)(x_0,\dots,x_n) = 0\cr
0 \leq i \leq n\end{matrix}\right.
\end{align}
has only the trivial solution $(x_0, \dots , x_n) = (0,\dots,0)$ in
$\C^{n+1}$. We remark that in this case this is true for the generic
$z \in \C^m$.

In order to prove our result for (weakly) general position (under the stronger assumption of pointwise general position this can be avoided), we finally will need some classical results on resultants,
see Lang \cite{b10}, section IX.3, for the precise definition, the existence and for the principal properties of resultants, as well as Eremenko-Sodin \cite{b4}, page 127:
Let $\big\{Q_j\big\}_{j=0}^n$ be a set of homogeneous
polynomials  of common degree $d \geq 1$ in
$\Cal K_f[x_0,\dots,x_n]$
 \begin{align*}
Q_j = \sum_{I \in \Cal T_d} a_{jI}x^I,\quad a_{jI} \in \Cal K_f \quad
(j = 0,\dots,n).
\end{align*}
Let $T = (\dots,t_{kI},\dots)$ \ ($k \in \{0,\dots,n\}$, $I \in \Cal T_d$)
be a family of variables. Set
\begin{align*}
\widetilde Q_j = \sum_{I \in \Cal T_d} t_{jI}x^I \in \Z[T,x],\quad
j = 0,\dots, n.
\end{align*}
Let $\widetilde R \in \Z[T]$ be the resultant of $\widetilde Q_0, \dots,
\widetilde Q_n$. This is a polynomial in the variables 
$T = (\dots,t_{kI},\dots)$ \ ($k \in \{0,\dots,n\}$, $I \in \Cal T_d$)
with integer coefficients, such that the condition 
$\widetilde R (T) =0$ is necessary and sufficient for the
existence of a nontrivial solution 
$(x_0, \dots , x_n) \not= (0,\dots,0)$ in $\C^{n+1}$
of the system of equations
\begin{align} \label{z}
\left\{ \begin{matrix}
\widetilde Q_{j}(T)(x_0,\dots,x_n) = 0\cr
0 \leq i \leq n\end{matrix}\right. \:.
\end{align}
From equations (\ref{z}) and (\ref{zz}) is follows immediately
that if $$\big\{Q_j= \widetilde Q_j(a_{jI})(x_0, \dots, x_n)\,, \:j=0, \dots , n\big\}$$ is an admissible set,
\begin{equation}R := \widetilde R(\dots, a_{kI}, \dots) \not\equiv 0\,.\label{zzz}
\end{equation}
Furthermore, since  $a_{kI} \in \Cal K_f$, we have  
$R \in \Cal K_f$.
 We finally will use the following result on resultants,
 which is contained in Theorem 3.4 in \cite{b10} (see also Eremenko-Sodin \cite{b4}, page 127, for a similar result):
\begin{proposition}\label{lang} There exists a positive integer $s$
and polynomials $\big\{\widetilde
b_{ij}\big\}_{0 \leq i, j \leq n}$ in $\Z[T,x]$, which are (without loss of generality) zero or
homogenous in $x$ of degree $s-d$,
such that
\begin{align*}
x_i^s \cdot \widetilde R = \sum_{j=0}^n \widetilde b_{ij} \widetilde Q_j\quad
\text{for all}\ i \in \{0,\dots,n\}.
\end{align*}
  \end{proposition}
\noindent If we still set
\begin{align*}
b_{ij} = \widetilde b_{ij}\big((\dots,a_{kI},\dots), (f_0,\dots,f_n)\big),\quad
0 \leq i, j \leq n,
\end{align*}
we get

\begin{align}\label{res}
f_i^s \cdot R = \sum_{j=0}^n b_{ij} \cdot Q_j(f_0,\dots,f_n)\quad
\text{for all}\ i \in \{0,\dots,n\}.
\end{align}
In particular, if $D \subset \C^m$ is a divisor contained in all divisors $f^{-1}(Q_j)$, $j=0,...,n$, then $R$ vanishes on $D$: This follows from (\ref{res}) since  $f=(f_0:...:f_n)$ is a reduced representation (and it follows in principle already directly from the definition of
the resultant).

\section{Main result}

Let $f,g$ be  nonconstant meromorphic maps of $\C^m$ into $\C P^n$.
Let $\big\{Q_j\big\}_{j=1}^q$ be an admissible set of homogeneous
polynomials in $\Cal K_f [x_0,\dots,x_n]$ with $\deg Q_j = d_j \geq
1$. Denote by $d, d^*,\tilde{d}$ respectively the
  least common multiple, the maximum number and the minimum number  of the $d_j$'s. Put $N=d\cdot(4(n+1)(2^n-1)(nd+1)+n+1)$. Set
$t_{\{Q_j\}_{j=1}^q}=1$ if the field $\Cal K_{\{Q_j\}_{j=1}^q}$
coincides with the complex number field $\C$ (ie. all $Q_j$ are fixed
hypersurface targets) and
$$
 t_{\{Q_j\}_{j=1}^q}
=\Bigg(\binom{n+N}{n}^2.\binom{q}{n}+\big[\frac{\big(
\binom{n+N}{n}^2.\binom{q}{n}-1 \big).\log\big(
\binom{n+N}{n}^2.\binom{q}{n}\big)}{\log(1+\frac{1}{4\binom{n+N}{n}N})}+1\big]^2\Bigg)^{
\binom{n+N}{n}^2.\binom{q}{n}-1}$$ if $\Cal
K_{\{Q_j\}_{j=1}^q}\ne\C,$
where we
denote $[x]:=\max\{k\in \Z: k\leq x\}$ for a real number $x.$
  Let
   $L=[\frac{d^*\cdot\binom{n+N}{n}t_{\{Q_j\}_{j=1}^q}-d^*}{d}+1].$

  With these notations, we state our main result:

 \begin{theorem}  \label {3.1} a) Assume that $f,g$ are algebraically nondegenerate
over $\Cal K_{\{Q_j\}_{j=1}^q}$ and satisfy \\
\indent i)  $\mathcal D^\alpha\big(\frac{f_k}{f_s}\big)=\mathcal
 D^\alpha\big(\frac{g_k}{g_s}\big)$ on
% $\big(\cup_{i=1}^qf^{-1}(Q_i)\big)\backslash \big(Zero (f_s\cdot g_s)\big),$ 
$\big(\cup_{i=1}^q(f^{-1}(Q_i)\cup g^{-1}(Q_i))\big)\backslash \big(Zero
(f_s.g_s)\big),$
 for  all  $|\alpha | < p,\ 0\leq k\ne s
 \leq n,$ where $p$ is a positive integer and $(f_0:\cdots:f_n),$ $(g_0:\cdots:g_n)$ are
 reduced representations of $f,g$ respectively.\\
Then for $q>n+\frac{2nL}{p\tilde{d}}+\frac{3}{2}$ , we have $f\equiv g.$\\
b) Assume that $f,g$ as in a) satisfy i) and\\
\indent ii) $\dim \big(f^{-1}(Q_i)\cap f^{-1}(Q_j)\big)\leq m-2$ for all
$1\leq i<j\leq q$.\\
Then for $q>n+\frac{2L}{p\tilde{d}}+\frac{3}{2}$ , we have $f\equiv g.$
\end{theorem}

We note that if $p=1$ the condition $i)$ becomes the following
usual condition: $f=g$ on $\cup_{i=1}^q(f^{-1}(Q_i)\cup g^{-1}(Q_i)),$ and we state this case again explicitly because of its importance:

 \begin{corollary}  \label {3.2} a) Assume that $f,g$ are algebraically nondegenerate
over $\Cal K_{\{Q_j\}_{j=1}^q}$ and satisfy\\ 
\indent i)  $f=g$ on $\cup_{i=1}^q(f^{-1}(Q_i)\cup g^{-1}(Q_i))$.\\
Then for $q>n+\frac{2nL}{\tilde{d}}+\frac{3}{2}$ , we have $f\equiv g.$\\
b) Assume that $f,g$ as in a) satisfy i) and\\
\indent ii) $\dim \big(f^{-1}(Q_i)\cap f^{-1}(Q_j)\big)\leq m-2$ for all
$1\leq i<j\leq q$.\\
Then for $q>n+\frac{2L}{\tilde{d}}+\frac{3}{2}$ , we have $f\equiv g.$
\end{corollary}

 In order to prove  Theorem 3.1, we need the following two results.
 The first one is similar to Lemma 5.1 in Ji \cite{J}, the second one
 is a special case of our main result in \cite{DT1}.
 \begin{proposition}  \label {3.3}Let $A_1,\dots,A_k$ be pure $(m-1)$- dimensional
 analytic subsets of $\C^m.$
 Let $f_1, f_2$ be  meromorphic  maps
of $\C^m$ into $\C P^n$. Then there exists a dense subset $\mathcal
C\subset \C^{n+1}\backslash \{0\}$ such that for any
$c=(c_0,\dots,c_n)\in \mathcal C$ the hyperplane $H_c$ defined by
$c_0w_0+\cdots+c_nw_n=0$ satisfies: $\dim\big(\cup_{j=1}^kA_j\cap
f_i^{-1}(H_c)\big)\leq m-2, i\in\{1,2\}.$
 \end{proposition}
\noindent {\bf Proof of Proposition~\ref{3.3}:}
 For any irreducible pure $(m-1)-$dimensional component $\sigma$ of
 $\cup_{j=1}^k A_j$  we set
\begin{equation*}
K_\sigma^i =\big\{(t_0,\dots,t_n) \in \C^{n+1} :
\sum\limits_{s=0}^nt_sf_{is} =0 \text{\ \ on\ }\sigma \big\}\ ,\quad
i\in\{1, 2\},
\end{equation*}
where $(f_{i0}:\cdots:f_{in})$ are reduced representations of
$f_i$. Then $K_\sigma^i$ is a complex vector subspace of $\C^{n+1}$.
Since $\text{dim}\{ f_{i0} = \cdots = f_{in} = 0\}\leq m-2,$ we get
that $\sigma\backslash\bigcup\limits_{i\in\{1, 2\}}\{ f_{i0} =
\cdots = f_{in} = 0\}\ne\varnothing.$ This implies that
   $\dim K_\sigma^i \leqslant n$. Let $K = \bigcup\limits_{i\in\{1, 2\}}\bigcup\limits_\sigma K_\sigma^i,$
   then $K$ is a union of at most a countable number of at 
   most $n-$dimensional complex vector subspaces in $\C^{n+1}$.
   Let $\mathcal{C} = \C^{n+1}\backslash K$.
    Then $\mathcal{C}$ meets the requirement of the Proposition.
 \qed

\begin{theorem}  \label {3.4} Under the same assumption as in Theorem~\ref{3.1},
we have
\begin{align*}
 (q-n-\frac{3}{2}) T_f(r) \leq \sum_{j=1}^q \frac{1}{d_j}
N^{(L)}_f(r,Q_j),
\end{align*}
for all $r \in [1, +\infty)$ excluding a Borel subset $E$ of $[1,
+\infty)$ with $\displaystyle{\int\limits_E} dr < +  \infty$.
\end{theorem}

\noindent {\bf Proof of Theorem~\ref{3.4}:} This is the special case of the Main Theorem and Proposition 1.2. in  \cite{DT1}
for $\epsilon = \frac{1}{2}$ and where we estimate the different
$d_j$'s in the numerators of the expressions entering into the
truncation level $L$ by $d^*$.
\qed
\\

\noindent {\bf Proof of Theorem~\ref{3.1}:} Assume that $f\not\equiv g.$ We first prove the following\\
{\bf Claim:} There
exist (fixed) hyperplanes $H_i: a_{i0}w_0+\dots a_{in}w_n=0\;(i=1,2)$ in $\C
P^n$ such that $S=S_{H_1, H_2}(f,g):=
\frac{H_1(f)}{H_2(f)}-\frac{H_1(g)}{H_2(g)}\not\equiv 0$ and
\begin{align}\dim(f^{-1}(Q_j)\cap f^{-1}\big(H_i)\big)\leq m-2,\; \dim(g^{-1}(Q_j)\cap g^{-1}\big(H_i)\big)\leq
m-2\label{1}
\end{align}
for all $j\in\{1,\dots,q\}$, $i\in\{1,2\}.$\\
{\bf Proof of the Claim}: By assumption i) of Theorem~\ref{3.1} we have pure $(m-1)-$dimensional analytic sets
\begin{align}A_j:=f^{-1}(Q_j)=g^{-1}(Q_j) \subset
\C^m, \:
j=1,\dots, q\,.\label{2}\end{align}
By Proposition~\ref{3.3}
 there exists a dense subset $\mathcal
C\subset \C^{n+1}\backslash \{0\}$ such that for any
$c=(c_0,\dots,c_n)\in \mathcal C$ the hyperplane $H_c$ defined by
$c_0w_0+\cdots+c_nw_n=0$ satisfies (\ref{1}), that is
$$\dim(A_j\cap f^{-1}\big(H_c)\big)\leq m-2,\; \dim(A_j\cap g^{-1}\big(H_c)\big)\leq
m-2$$
for all $j\in\{1,\dots,q\}$. Since $f,g$ are algebraically nondegenerate
over $\Cal K_{\{Q_j\}_{j=1}^q}$, so in particular algebraically nondegenerate over $\C$, we have that $L_c(f)\not\equiv 0$
and $L_c(g)\not\equiv 0$ are holomorphic functions for all 
$c=(c_0,\dots,c_n)\in \mathcal C$, where 
$L_c(f):= \sum_{i=0}^nc_if_i$ with a reduced representation
$f=(f_0:\dots :f_n)$ and 
$L_c(g):= \sum_{i=0}^nc_ig_i$ with a reduced representation
$g=(g_0:\dots :g_n)$. Finally for 
$c^{(1)}, c^{(2)}\in \mathcal C$, we put 
 $S_{c^{(1)},c^{(2)}}(f,g):=
\frac{L_{c^{(1)}}(f)}{L_{c^{(2)}}(f)}-\frac{L_{c^{(1)}}(g)}{L_{c^{(2)}}(g)}$. In order to prove the Claim it suffices to show that for some $c^{(1)}, c^{(2)}\in \mathcal C$, $S_{c^{(1)},c^{(2)}}(f,g) \not\equiv 0$. Assume the contrary. Then for all $0 \leq i < j \leq n$
there exist sequences $(c^{(1)})_{\nu}$, $(c^{(2)})_{\nu}$, 
$\nu \in \N$, of elements in $\mathcal C$ such that 
$L_{(c^{(1)})_{\nu}}(f) \rightarrow f_i$ and
$L_{(c^{(2)})_{\nu}}(f) \rightarrow f_j$.
From this we get $$0\equiv
S_{(c^{(1)})_{\nu},(c^{(2)})_{\nu}}(f,g) \rightarrow \frac{f_i}{f_j} -
\frac{g_i}{g_j}\,,$$
what implies $0 \equiv \frac{f_i}{f_j} -
\frac{g_i}{g_j}$ for all $0\leq i<j\leq n$, contradicting the assumption $f\not\equiv g$. This proves the claim.\qed

Since $f=g$ on $\cup_{j=1}^q f^{-1}(Q_j),$ for any generic point
$$z_0\in\cup_{j=1}^q
f^{-1}(Q_j)\backslash \big( f^{-1}(H_2)\cup g^{-1}(H_2)\big)$$ 
(outside an analytic subset of codimension at least 2), there exists
$s\in\{0,\dots,n\}$ such that both of $f_s(z_0)$ and $g_s(z_0)$ are
different from zero. Then by assumption i) we have
\begin {align*}
\mathcal D^\alpha S(z_0)&=\mathcal D^\alpha\bigl (
 \frac{H_1(f)}{H_2(f)}-\frac{H_1(g)}{H_2(g)}\bigl )(z_0)\\&
=\mathcal D^\alpha\bigl (
 \frac{\sum_{v=0}^{n}\frac{f_v}{f_s}a_{1v}}{\sum_{v=0}^{n}\frac{f_v}{f_s}a_{2v}}-\frac{\sum_{v=0}^{n}\frac{g_v}{g_s}a_{1v}}{\sum_{v=0}^{n}\frac{g_v}{g_s}a_{2v}}\bigl
 )(z_0)=0
\end{align*}
for all $ |\alpha |< p.$

\noindent This implies that

\begin{align}
\nu_S\ge p \;\text{on}\;\cup_{j=1}^q f^{-1}(Q_j)\backslash
\big(A\cup f^{-1}(H_2)\cup g^{-1}(H_2)\big).\label{3}
\end{align}
where $A$ is an analytic subset of codimension at least 2.

Now we will estimate the divisors $\nu_{Q_j \circ f}$ by making use of the resultants:
In fact, for any $J=\{j_0,...,j_n\}  \subset \{1,2,...,q\}$, let $R_J$ be the resultant of
of $Q_{j_0},...,Q_{j_n}$. Then if $D \subset \C^m$ is a divisor contained in all divisors $f^{-1}(Q_{j_k})$, $k=0,...,n$, then $R_J$ vanishes on $D$. Thus, we get
\begin{align}\label{n}
 \sum_{j=1}^q\text{min}\{1, \nu_{Q_j \circ f}\} \leq n \cdot \text{min}\{1, \sum_{j=1}^q \nu_{Q_j \circ f}\} + (q-n) \cdot \text{min}\{1,  \sum_{|J|=n+1}\nu_{R_J}\}
\end{align}

By (\ref{1}), (\ref{2}),(\ref{3}),  (\ref{n}), by the First Main Theorem and since $R_J \in {\cal K}_f$, we have
\begin{align}
\sum_{j=1}^qN_g^{(1)}(r,Q_j)=\sum_{j=1}^qN_f^{(1)}(r,Q_j) &\leq
\frac{n}{p}N_S(r) + o(T_f(r))\label{4a}
\end{align}
Furthermore, by the First Main Theorem
\begin{align}
N_S(r)
&\leq
T_{\frac{H_1(f)}{H_2(f)}-\frac{H_1(g)}{H_2(g)}}(r)+O(1)\notag\\
&\leq
T_{\frac{H_1(f)}{H_2(f)}}(r)+T_{\frac{H_1(g)}{H_2(g)}}(r)+O(1)\notag\\
&\leq T_f(r)+T_g(r)+O(1).\label{4}
\end{align}
Thus,
\begin{align}
 \sum_{j=1}^q\big(N_f^{(1)}(r,Q_j)+N_g^{(1)}(r,Q_j)\big)\leq
 \frac{2n}{p}\big(T_f(r)+T_g(r)\big)+o(T_f(r)).\label{5}
\end{align}
By Theorem~\ref{3.4}  and by the First Main Theorem, we have
\begin{align}
 (q-n-\frac{3}{2})T_f(r)\leq\sum_{j=1}^q\frac{1}{d_j}N_f^{(L)}(r,Q_j)\notag\\
 \leq\sum_{j=1}^q\frac{L}{d_j}N_f^{(1)}(r,Q_j)=\sum_{j=1}^q\frac{L}{d_j}N_g^{(1)}(r,Q_j)\leq
 qLT_g(r)+o(T_f(r))\label{6}
\end{align}
for all $r \in [1, +\infty)$ excluding a Borel subset $E$ of $(1,
+\infty)$ with $\displaystyle{\int\limits_E} dr < +  \infty$
(note that $Q_j\in\mathcal K_f[x_0,\dots,x_n]).$

\noindent This implies that $\mathcal K_f\subset\mathcal K_g.$ Then
$\{Q_j\}_{j=1}^q\subset\mathcal K_g[x_0,\dots,x_n].$
 So we can
apply Theorem~\ref{3.4} for both  meromorphic maps $f$ and $g$ with moving
hypersurfaces $\{Q_j\}_{j=1}^q.$ By Theorem~\ref{3.4} and by the First Main
Theorem, we have
\begin{align}
 (q-n-\frac{3}{2})\big(T_f(r)+T_g(r)\big)\leq\sum_{j=1}^q\frac{1}{d_j}\big(N_f^{(L)}(r,Q_j)+N_g^{(L)}(r,Q_j)\big)\notag\\
 \leq\frac{L}{\tilde{d}}\sum_{j=1}^q\big(N_f^{(1)}(r,Q_j)+N_g^{(1)}(r,Q_j)\big)\label{7}
 \end{align}
for all $r \in [1, +\infty)$ excluding a Borel subset $E$ of $(1,
+\infty)$ with $\displaystyle{\int\limits_E} dr < +  \infty$.

Combining with (\ref{5}), we get
\begin{align}
 (q-n-\frac{3}{2})\big(T_f(r)+T_g(r)\big)\leq\frac{2nL}{p\tilde{d}}\big(T_f(r)+T_g(r)\big)+o(T_f(r))
\label{8}\end{align}
for all $r \in [1, +\infty)$ excluding a Borel subset $E$ of $(1,
+\infty)$ with $\displaystyle{\int\limits_E} dr < +  \infty$.
This is a contradiction, since
$q>n+\frac{2nL}{p\tilde{d}}+\frac{3}{2}$, thus finishing the proof of part a).\\

In order to prove b), we observe that under the additional assumption  ii), we can
improve (\ref{4a}), namely we get, by using (\ref{1}), (\ref{2}), (\ref{3}) and assumption ii)
\begin{align}
\sum_{j=1}^qN_g^{(1)}(r,Q_j)=\sum_{j=1}^qN_f^{(1)}(r,Q_j) &\leq
\frac{1}{p}N_S(r)\label{4ab}
\end{align}
This improves (\ref{5}), namely we get from (\ref{4}) and (\ref{4ab}):
\begin{align}
 \sum_{j=1}^q\big(N_f^{(1)}(r,Q_j)+N_g^{(1)}(r,Q_j)\big)\leq
 \frac{2}{p}\big(T_f(r)+T_g(r)\big)+O(1).\label{5b}
\end{align}
Using this (\ref{8}) becomes, by using now (\ref{7}) and (\ref{5b}):
\begin{align}
 (q-n-\frac{3}{2})\big(T_f(r)+T_g(r)\big)\leq\frac{2L}{p\tilde{d}}\big(T_f(r)+T_g(r)\big)+O(1)
\label{8b}\end{align}
for all $r \in [1, +\infty)$ excluding a Borel subset $E$ of $(1,
+\infty)$ with $\displaystyle{\int\limits_E} dr < +  \infty$.
This is a contradiction, since
$q>n+\frac{2L}{p\tilde{d}}+\frac{3}{2}$, thus finishing the proof of part b).
\qed

\vspace{1cm}

 \noindent  Gerd Dethloff$^{1-2} $ \\
 $^1$ Universit\'e Europ\'eenne de Bretagne, France\\
 $^2$
Universit\'{e} de Brest \\
   Laboratoire de math\'{e}matiques \\
UMR CNRS 6205\\
6, avenue Le Gorgeu, BP 452 \\
   29275 Brest Cedex, France \\
e-mail: gerd.dethloff@univ-brest.fr\\

\noindent Tran Van Tan\\
Department of Mathematics\\
  Hanoi National University of Education\\
 136-Xuan Thuy street, Cau Giay, Hanoi, Vietnam\\
e-mail: tranvantanhn@yahoo.com; vtran@ictp.it

\end{document}